\documentclass{article}
\usepackage{amsthm}
\usepackage{amsmath}
\usepackage{amsfonts}
\usepackage[letterpaper,body={14.5cm,23cm}, mag=1000]{geometry}
\usepackage{amssymb}
\usepackage{secdot}
\usepackage{setspace}
\usepackage{titlesec}
\titleformat{\section}[runin]{\bfseries\filcenter}{\thesection}{1em}{}

\renewcommand{\thesection}{\arabic{section}}
\title{\large \bf A note on commuting automorphisms of some finite $p$-groups}
\author{\small \bf Sandeep Singh\footnote{Research supported by CSIR, New Delhi} and Deepak Gumber\footnote{Research supported by NBHM, Department of Atomic Energy, Govt. of India}\\
\small \em School of Mathematics and Computer Applications\\
\small \em Thapar University, Patiala - 147 004,
India\\
} 
\date{}

\DeclareMathOperator{\Aut}{Aut}

\newtheorem{thm}{Theorem}

\begin{document}
\maketitle

\begin{abstract}
\noindent An automorphism $\alpha$ of a group $G$ is called a commuting automorphism if each element $x$ in $G$ commutes with its image $\alpha(x)$ under $\alpha$. Let $A(G)$ denote the set of all commuting automorphisms of $G$. Rai [Proc. Japan Acad., Ser. A {\bf 91} (2015), no. 5, 57-60] has given some sufficient conditions on a finite $p$-group $G$ such that $A(G)$ is a subgroup of $\Aut(G)$ and, as a consequence, has proved that in a finite $p$-group $G$ of co-class 2, where $p$ is an odd prime, $A(G)$ is a subgroup of $\Aut(G)$. We give here very elementary and short proofs of main results of Rai.
\end{abstract}

\vspace{2ex}

\vspace{2ex}


\vspace{2ex}

Let $G$ be a finite non-abelian $p$-group of order $p^n$, and let $\gamma_k(G)$ and $Z_k(G)$ respectively denote the $k$th terms of the lower and upper central series of $G$.  For convenience, $\gamma_2(G)$ and $Z_1(G)$ are respectively denoted as $G'$ and $Z(G)$. We call $G$ an $A(G)$-group if $A(G)$ is a subgroup of $\Aut(G)$. Let $\alpha,\beta\in A(G)$. Then, for any $x\in G^2$, $x^{-1}\alpha(x),x^{-1}\beta(x)\in Z_2(G))$,  by \cite[Theorem 1.4]{dea}. Now if $p$ is an odd prime, then $G^2=G$ and hence $x^{-1}\alpha(x),x^{-1}\beta(x)\in Z_2(G)$. If $Z_2(G))$ is abelian, then $[\alpha(x),\beta(x)]=[x^{-1}\alpha(x),x^{-1}\beta(x)]=1$ and hence $A(G)$ is a subgroup of $\Aut(G)$ by \cite[Lemma 2.4(vi)]{dea}. Suppose that
$|Z_2(G)/Z(G)|=p^2$ and $Z(G)=\gamma_k(G)$ for some $k\ge 2$.
Then $G$ is of nilpotence class $k$ and hence $\gamma_{k-1}(G)\le Z_2(G)$. If $k=2$, then $G$ is an $A(G)$-group by \cite[Lemma 2.2]{vos}. Assume that $k\ge 3$. Since $\gamma_{k-1}(G)$ commutes with $Z_2(G)$, $\gamma_{k-1}(G)$ is a central subgroup of $Z_2(G)$. It follows that $Z_2(G)$ is abelian, because $Z_2(G)/\gamma_{k-1}(G)$ is cyclic. We have thus proved the following main theorem of Rai.
\begin{thm} [{\cite[Theorem 3.3]{rai}}]
Let $p$ be an odd prime and $G$ be a finite $p$-group such that $|Z_2(G)/Z(G)|=p^2$ and $Z(G)=\gamma_k(G)$ for some $k\ge 2$. Then $G$ is an $A(G)$-group.
\end{thm}
Observe that if $G$ is of co-class 2, then $|Z_2(G)/Z(G)|=p$ or $p^2$. It follows that $Z_2(G)$ is abelian as explained above. We thus have the following theorem of Rai.
\begin{thm} [{\cite[Theorem A]{rai}}]
Let $G$ be a finite $p$-group of co-class $2$ for an odd prime $p$. Then $G$ is an $A(G)$-group.
\end{thm}

\end{document}